\documentclass[12pt]{article}
\usepackage[centertags]{amsmath}
\usepackage{amsfonts}
\usepackage{amssymb}
\usepackage{amsthm}

\theoremstyle{plain}
\newtheorem{thm}{Theorem}[section]
\newtheorem{cor}[thm]{Corollary}
\newtheorem{lem}[thm]{Lemma}
\newtheorem{prop}[thm]{Proposition}
\newtheorem{remark}[thm]{Remark}
\theoremstyle{definition}
\newtheorem{defn}[thm]{Definition}
\theoremstyle{remark}
\numberwithin{equation}{section}

\def\sqr#1#2{{\vcenter{\hrule height.#2pt

       \hbox{\vrule width.#2pt height#1pt \kern#1pt

           \vrule width.#2pt}

       \hrule height.#2pt}}}

\def\openC{{\bf\rm C}\kern-.4em {\vrule height 1.4ex width .08em depth-.04ex}}

\newcommand{\beast}{\begin{eqnarray*}}
\newcommand{\eeast}{\end{eqnarray*}}

\title{ Distance-regular graphs, 
pseudo primitive idempotents, and  
the Terwilliger algebra}

\author{Paul Terwilliger
\and Chih-wen Weng
}

\date{}
\begin{document}
\maketitle

\begin{abstract}

\bigskip

Let $\Gamma$ denote a distance-regular graph with diameter $D\geq
3$, intersection numbers $a_i,b_i,c_i$ and Bose-Mesner algebra $ \mathbf{M}.$
For
 $\theta\in \mathbb{C} \cup
\infty$ we define a 1 dimensional subspace of 
 $\mathbf{M}$
which we call $\mathbf{M}(\theta)$.
If $\theta\in \mathbb{C}$ then
$\mathbf{M}(\theta)$ consists of those
$Y $ in $\mathbf{M} $ 
such that $(A-\theta I)Y\in
\mathbb{C} A_D$,
where $A$ (resp. $A_D$) is the adjacency matrix (resp. $D$th distance matrix)
of
$\Gamma.$ If $\theta = \infty$
then $\mathbf{M}(\theta)=
\mathbb{C}A_D.$
By a {\it pseudo primitive idempotent} for $\theta$ we mean a
nonzero element of $\mathbf{M}(\theta).$
We use these
as follows.
Let $X$ denote the vertex
set of $\Gamma$ and fix $x\in X.$ Let $ \mathbf{T}$
denote the subalgebra of ${\rm Mat}_X(\mathbb{C})$ generated by
$A,$ $E_0^*, E_1^*, \cdots , E_D^*,$ where
 $E_i^*$ denotes the projection
onto the $i$th subconstituent of $\Gamma$ with respect to
$x.$ $\mathbf{T}$ is called the Terwilliger algebra.
%
Let $W$ denote an irreducible 
$\mathbf{T}$-module.
By the {\it endpoint} of
$W$ we mean ${\rm min}\{i|E_i^*W\not=0\}.$
$W$ is called {\it thin} whenever
${\rm
dim}(E_i^*W)\leq 1$ for $0\leq i\leq D$.
Let $V
=\mathbb{C}^X$
 denote the standard 
$\mathbf{T}$-module.
Fix 
$0\not=v\in E_1^*V$ with $v$ orthogonal to the all 1's vector.
 We define
$(\mathbf{M}; v):=\{P\in \mathbf{M}| Pv\in E_D^* V\}$.
We show the following are equivalent: (i) 
$\hbox{dim}(\mathbf{M}; v)\geq 2$;
(ii) $v$ is contained in a
thin irreducible $ \mathbf{T}$-module with endpoint $1.$
Suppose (i), (ii) hold. We show
$(
\mathbf{M};v)$ 
has a basis $J,E$ where
$J$
 has all entries 1 and
$E$ is defined as follows.
Let $W$ denote the 
 $ \mathbf{T}$-module  which satisfies
(ii). Observe $E^*_1W$ is an eigenspace for $E^*_1AE^*_1$;
let $\eta$ denote the corresponding eigenvalue.
Define $\widetilde \eta = -1 -b_1(1+\eta)^{-1}$ if
$\eta \not=-1$ and 
$\widetilde \eta = \infty $ if $\eta = -1$. Then
$E$ is a pseudo primitive idempotent for $\widetilde \eta$.

\end{abstract}

\section{Introduction}\label{s1}
\bigskip

\noindent Let $\Gamma$ denote a distance-regular graph with
diameter $D\geq 3$, intersection numbers $a_i,b_i, c_i$,
Bose-Mesner algebra $\mathbf{M}$
 and path-length distance
function $\partial$ (see Section 2 for formal definitions). In
order to state our main theorems we make a few comments. Let $X$
denote the vertex set of $\Gamma.$  Let $V= \mathbb{C}^X$ denote
the vector space over $ \mathbb{C}$ consisting of column vectors
whose coordinates are indexed by $X$ and whose entries are in 
 $ \mathbb{C}$.
 We endow $V$ with the Hermitean inner product
$\langle\,,\,\rangle$
satisfying
 $ \langle u,v\rangle=u^t\overline v$ for all
$u,v\in V$. For each $y\in X$ let $\hat y$ denote the vector in
$V$ with a $1$ in the $y$ coordinate and $0$ in all other
coordinates. We observe $\{\hat y|y\in X\}$ is an orthonormal
basis for $V.$ Fix  $x \in X$. For $0 \leq i \leq D$
let $E^*_i$ denote the diagonal matrix in ${\rm
Mat}_X(\mathbb{C})$ which has $yy$ entry $1$ (resp. $0$) whenever
$\partial(x,y)=i$ (resp. $\partial(x,y)\not=i$). We observe
$E^*_i$ acts on $V$ as the projection onto the $i$th
subconstituent of $\Gamma$ with respect to $x$. For $0 \leq i \leq
D$ define $s_i =\sum  {\hat y}$, where the sum is over all
vertices $y \in X$ such that $\partial(x,y)=i$. We observe $s_i
\in E^*_iV$. Let $v$ denote a
nonzero
vector in $E^*_1V$
 which is orthogonal to $s_1$.
 We define
$$(\mathbf{M};v):=\{P\in \mathbf{M}\ |\ Pv\in E_D^*V\}.$$ We
observe $(\mathbf{M};v)$ is a subspace of $\mathbf{M}$.
\medskip
\noindent We consider the dimension of $(\mathbf{M};v)$. We first
observe $(\mathbf{M};v) \not=0$. To see this, let $J $ denote the
matrix in ${\rm Mat}_X(\mathbb{C})$ which has all entries 1. It is
known $J$ is contained in $\mathbf{M}$ \cite[p. 64]{BanIto}. In fact $J
\in (\mathbf{M};v)$; the reason is $Jv =0$ since $v$ is orthogonal
to $s_1$. Apparently $(\mathbf{M};v)$ is nonzero so it has
 dimension
at least 1. We now consider when
does $(\mathbf{M};v)$ have dimension
at least 2?
To answer this question we recall the Terwilliger algebra.
Let
$\mathbf{T}$
  denote the subalgebra of
 ${\rm Mat}_X(\mathbb{C})$
generated by $A, E_0^*, E_1^*,\ldots,E_D^*$, where $A$ denotes the
adjacency matrix of $\Gamma$. The algebra $\mathbf{T}$ is known as
the {\it Terwilliger} algebra (or  {\it subconstituent } algebra)
of $\Gamma$ with respect to $x$ 
\cite{te92,te93,te93-2}.
By a
$\mathbf{T}$-{\it module} we mean a subspace $W \subseteq V$ such
that $\mathbf{T}W \subseteq W$. Let $W$ denote a
$\mathbf{T}$-module. We say $W$ is {\it irreducible} whenever
$W\not=0$
and
$W$
does not contain a $\mathbf{T}$-module other than $0$ and $W$. Let
$W$ denote an irreducible $\mathbf{T}$-module.
 By the {\it endpoint} of $W$ we mean the minimal
integer $i$ $(0 \leq i \leq D)$ such that $E_i^*W\not=0$.
We say $W$ is {\it
thin} whenever $E^*_iW$ has dimension at most 1 for $0 \leq  i
\leq D$. 
 We now
state our main theorem.

\begin{thm}\label{10-20i}
Let $v$ denote a nonzero vector in $E^*_1V$ which is orthogonal to
$s_1$.
Then the following (i), (ii) are equivalent.
\begin{enumerate}
\item[(i)]
$(\mathbf{M}; v)$ has dimension at least $2.$
\item[(ii)] $v$ is contained in a thin irreducible $\mathbf{T}$-module
with endpoint 1.
\end{enumerate}
Suppose (i),(ii) hold above. Then $(\mathbf{M}; v)$ has dimension
exactly $2.$
\end{thm}
\noindent
With reference to Theorem~\ref{10-20i}, suppose for the moment
that (i), (ii) hold. We find a basis for $(\mathbf{M}; v)$. To
describe our basis we need some notation. Let $\theta_0 > \theta_1
> \cdots > \theta_D$ denote the distinct eigenvalues of $A$, and
for $0\leq i  \leq D$ let $E_i$ denote the primitive idempotent of
$\mathbf{M} $ associated with $\theta_i$. We recall $E_i$
satisfies $(A-\theta_iI)E_i=0$.
We introduce a
type of element in
$\mathbf{M} $ which generalizes the $E_0, E_1, \ldots,
E_D$. We call this type of element a {\it pseudo primitive
idempotent} for $\Gamma$.
In order to define
 the pseudo primitive idempotents,
we first define for each
 $\theta\in \mathbb{C} \cup
\infty$ a  subspace
of
 $\mathbf{M}$ which we call
$\mathbf{M}(\theta)$.
For $\theta\in \mathbb{C}$, $\mathbf{M}(\theta)$ consists of those
elements $Y $ of $\mathbf{M} $ such that $(A-\theta I)Y\in
\mathbb{C} A_D$, where $A_D$ is the $D$th distance matrix of
$\Gamma.$
We define $\mathbf{M}(\infty)=
\mathbb{C}A_D.$
We show
$\mathbf{M}(\theta)$ has dimension 1 for
all
 $\theta\in \mathbb{C} \cup
\infty$.
Given distinct
 $\theta, \theta'$ in $\mathbb{C} \cup
\infty$,
we show
$\mathbf{M}(\theta) \cap
\mathbf{M}(\theta')=0.
$
For $0 \leq i \leq D$ we show $\mathbf{M}(\theta_i) =
 \mathbb{C} E_i$.
Let $\theta\in \mathbb{C}\cup \infty.$ By a {\it
pseudo primitive idempotent} for $\theta,$ we mean a
nonzero element of $\mathbf{M}(\theta).$
%
Before proceeding we define an
 involution
 on
 $\mathbb{C}\cup \infty $.
For $\eta
\in
 \mathbb{C}\cup \infty $ we define
\beast \widetilde{\eta}=
 \begin{cases}
    \infty\ \ \  & \text{if $\eta=-1$},\\
    -1\ \ \      & \text{if $\eta=\infty$},\\
    -1-\frac{b_1}{1+\eta}\ \ \  & \text{if $\eta \not=-1, \eta \not=\infty$}.
  \end{cases}
\eeast We observe $\widetilde{\widetilde{\eta}}=\eta$\ for
$\eta\in\mathbb{C}\cup \infty$.
\medskip
\noindent
 Let $W$ denote a thin irreducible
$\mathbf{T}$-module with endpoint $1.$ Observe $E^*_1W$ is a one
dimensional eigenspace for $E^*_1AE^*_1$; let $\eta$ denote the
corresponding eigenvalue. We call $\eta$ the {\it local
eigenvalue} of $W$.
\begin{thm}\label{11-23}
Let $v$ denote a nonzero vector in $E^*_1V$ which is
orthogonal to $s_1$. Suppose $v$ satisfies the equivalent
conditions (i), (ii) in
Theorem~\ref{10-20i}.
Let $W$ denote the
 $\mathbf{T}$-module
from
part (ii) of that theorem
 and let
$\eta $ denote the local eigenvalue for $W$. Let $E$ denote a
pseudo primitive idempotent for $\widetilde \eta$. Then $J, E$
form a basis for $(\mathbf{M}; v)$.
\end{thm}
\bigskip

\noindent
We comment on when the scalar $\widetilde \eta $ 
from  
 Theorem 
\ref{11-23}
 is an eigenvalue of $\Gamma$. Let $W$ denote a
thin irreducible
 $\mathbf{T}$-module with endpoint 1 and local eigenvalue $\eta$.
It is known ${\widetilde \theta_1} \leq \eta \leq {\widetilde
\theta_D}$ \cite[Theorem 1]{te86}. If $\eta= {\widetilde \theta_1}$ then
 $\widetilde \eta =\theta_1$.
If $\eta= {\widetilde \theta_D}$ then
 $\widetilde \eta =\theta_D$.
We show that if  ${\widetilde \theta_1} < \eta < {\widetilde \theta_D}$ then
$\widetilde \eta$ is not an eigenvalue of $\Gamma$.
\bigskip

\noindent The paper is organized as follows. In section~\ref{s2}
we give some
preliminaries on distance-regular graphs. In section~\ref{s4}
and section~\ref{s5} we review some basic results
on the Terwilliger algebra and its modules. We prove
Theorem~\ref{10-20i} in section~\ref{main1}.
In
 section~\ref{s3} we discuss
pseudo primitive idempotents. In
 section~\ref{s6}
we discuss local
eigenvalues.
 We prove
Theorem~\ref{11-23} in section~\ref{main2}.
\section{Preliminaries}\label{s2}

\noindent In this section we review some definitions and basic
concepts. See the books by Bannai and Ito \cite{BanIto} or Brouwer,
Cohen, and Neumaier \cite{bcn} for more  background information.
\bigskip

\noindent
Let $X$ denote a nonempty  finite  set.
Let
 $\hbox{Mat}_X({\mathbb C})$ 
denote the $\mathbb  C$-algebra
consisting of all matrices whose rows and columns are indexed by $X$
and whose entries are in $\mathbb C  $. Let
$V={\mathbb C}^X$ denote the vector space over $\mathbb C$
consisting of column vectors whose 
coordinates are indexed by $X$ and whose entries are
in $\mathbb C$.
We observe
$\hbox{Mat}_X({\mathbb C})$ 
acts on $V$ by left multiplication.
We endow $V$ with the Hermitean inner product $\langle \, , \, \rangle$ 
which satisfies
$\langle u,v \rangle = u^t\overline{v}$ for all $u,v \in V$,
where $t$ denotes transpose and $-$ denotes complex conjugation.
For all $y \in X,$ let $\hat{y}$ denote the element
of $V$ with a 1 in the $y$ coordinate and 0 in all other coordinates.
We observe $\{\hat{y}\;|\;y \in X\}$ is an orthonormal basis for $V.$

\medskip
\noindent 
Let $\Gamma=(X,R)$ denote a finite, undirected, connected graph
without loops or multiple edges, with vertex set $X$, edge set
$R$, path-length distance function $\partial$ and diameter
$D:={\rm max}\{\partial(x,y)|x,y\in X\}.$ We say $\Gamma$ is {\it
distance-regular} whenever for all integers $h,i,j$ $(0\leq
h,i,j\leq D)$ and for all $x,y\in X$ with $\partial (x,y)=h,$ the
number
\begin{equation}
p^h_{ij}=|\{z\in X|\partial (x,z)=i, \partial (z,y)=j\}|\label{1}
\end{equation}
is independent of $x$ and $y.$ The integers $p^h_{ij}$ are called
the {\it intersection numbers} for $\Gamma$. Observe
$p^h_{ij}=p^h_{ji}$ $(0\leq h,i,j\leq D).$ We abbreviate
$c_i:=p^i_{1 i-1}$ $(1\leq i\leq D)$, $a_i:=p^i_{1i}$ $(0\leq
i\leq D)$, $b_i:=p^i_{1 i+1}$ $(0\leq i\leq D-1)$, $k_i:=p^0_{ii}$
$(0\leq i\leq D)$, and for convenience we set $c_0:=0$ and
$b_D:=0$. Note that $b_{i-1}c_i\not=0$ $(1\leq i\leq D).$
\bigskip

For the rest of this paper we assume $\Gamma=(X,R)$ is
distance-regular with diameter $D\geq 3$. By (\ref{1}) and the
triangle inequality,
\begin{eqnarray}\label{2}
p^h_{i1}&=&0\ \ \ \ \ {\rm if}\ \ |h-i|>1\ \ \ \ \ (0\leq h,i\leq
D),
\\
p^1_{ij}&=&0\ \ \ \ \ {\rm if}\ \ |i-j|>1\ \ \ \ \ (0\leq i,j\leq
D).
\label{3_30.3}
\end{eqnarray}
Observe $\Gamma$ is regular with valency $k=k_1=b_0$, and that
$k=c_i+a_i+b_i$ for $0\leq i\leq D$.
By \cite[p. 127]{bcn} we have
\begin{eqnarray}\label{4}
k_{i-1}b_{i-1}=k_ic_i\ \ \ \ \ (1\leq i\leq D).
\end{eqnarray}
\noindent 
We recall the Bose-Mesner algebra of $\Gamma$.
For $0\leq i 
\leq D$ let $A_i$ denote the matrix in $\hbox{Mat}_X(\mathbb{C})$
which has 
$yz$ entry
\begin{equation}
\nonumber
(A_i)_{yz}=\begin{cases}
    1  & \text{if $\partial (y,z)=i$}\\
    0  & \text{if $\partial (y,z)\not=i$}
    \end{cases}\ \ \ \ \ (y,z\in X).
\end{equation}
We call $A_i$ the {\it $i$th distance matrix} of $\Gamma.$
For notational convenience we define $A_i=0$ for $i<0$ and $i>D$.
Observe
(ai) $A_0=I$; (aii) $\sum_{i=0}^D A_i=J$;
(aiii) ${\overline{A_i}}=A_i$  $(0\leq i\leq D)$;
(aiv) $A^t_i=A_i$ $(0\leq i\leq D)$,
(av) $A_iA_j=\sum_{h=0}^D p^h_{ij}A_h$ $(0\leq i,j \leq D)$,
where $I$ denotes the identity matrix and
$J$ denotes the all ones matrix.
We abbreviate $A:=A_1$ and call this
the {\it adjacency matrix} of $\Gamma$.
%
Let $\mathbf{M}$ denote the subalgebra of ${\rm
Mat}_X(\mathbb{C})$ generated by $A$.
Using (ai)--(av)  we find
$A_0, A_1, \cdots, A_D$ form a basis of $\mathbf{M}$.
We call $\mathbf{M}$ 
the {\it Bose-Mesner algebra} of $\Gamma.$ 
By
\cite[p. 59, p. 64]{BanIto}, $\mathbf{M}$ has a second basis $E_0$,
$E_1$, $\cdots$, $E_D$ such that
(ei) $E_0=|X|^{-1}J$;
(eii) $\sum_{i=0}^D E_i=I$; (eiii) ${\overline{E_i}}=E_i
$ $(0\leq i\leq D)$;
(eiv) $E_i^t=E_i$ $(0\leq 
i\leq D)$;
(ev)
$E_iE_j=\delta_{ij}E_i $ $(0\leq
i,j\leq D)$.
We call $E_0,$ $E_1,$ $\cdots$, $E_D$  the {\it primitive
idempotents} for $\Gamma$.
 Since $E_0,$ $E_1,$ $\cdots$, $E_D$
form a
basis for 
 $\mathbf{M}$
 there exists complex scalars
$\theta_0,\theta_1,\cdots,\theta_D$ such that
  $A=\sum_{i=0}^D\theta_iE_i$.
By this and 
(ev) we find
$AE_i=\theta_i E_i$ for $0\leq i\leq D$.
Using (aiii) and (eiii) we find each of
 $\theta_0,$ $\theta_1,$ $\cdots$, $\theta_D$ 
is a real number. 
 Observe $\theta_0$, $\theta_1$, $\cdots$, $\theta_D$ are
mutually distinct since $A$ generates $\mathbf{M}.$ By
\cite[p.197]{BanIto} we have  $\theta_0=k$ and $-k\leq \theta_i\leq k$
for 
 $0\leq i\leq D.$  Throughout this paper, 
we assume
$E_0,$ $E_1,$ $\cdots$, $E_D$ 
are indexed so that
 $\theta_0>\theta_1>\cdots>\theta_D.$
We call $\theta_i$ the {\it $i$th eigenvalue} of $\Gamma$.

\medskip
\noindent We recall some polynomials. To motivate 
these 
we make a comment.
Setting $i=1$ in (av) and using (\ref{2}),
\begin{equation}\label{10}
  AA_j=b_{j-1}A_{j-1}+a_jA_j+c_{j+1}A_{j+1}\ \ \ \ \ (0\leq j\leq
  D-1),
\end{equation}
where $b_{-1}=0$.
Let $\lambda$ denote an indeterminate and let
 $\mathbb{C}[\lambda]$ denote the $\mathbb{C}$-algebra
consisting of  all polynomials in $\lambda$ which have coefficients in
$\mathbb{C}.$  Let $f_0,$ $f_1,$ $\cdots$, $f_D$ denote the
polynomials in $\mathbb{C}[\lambda]$  which satisfy $f_0=1$ and
\begin{equation}\label{14}
  \lambda f_j=b_{j-1}f_{j-1}+a_jf_j+c_{j+1}f_{j+1}\ \ \ \ \
  (0\leq j\leq D-1),
\end{equation}
where $f_{-1}=0$.
For $0 \leq j \leq D$ the degree of
 $f_j$ is  exactly $j$.
  Comparing 
  (\ref{10})
  and
  (\ref{14})
  we find
 $A_j=f_j(A)$.

\section{The Terwilliger algebra}\label{s4}
For the remainder of this paper we fix $x\in X$. For
$0\leq i\leq D$ let $E_i^*=E_i^*(x)$ denote the
diagonal matrix in ${\rm Mat}_X( \mathbb{C})$ which has
$yy$ entry
\begin{equation}\label{29}
(E_i^*)_{yy}=
  \begin{cases}
    1 & \text{if $\partial (x,y)=i$} \\
    0 & \text{if $\partial (x,y)\not=i$}
  \end{cases}
\ \ \ \ \ (y\in X).
\end{equation}
We call $E_i^*$ the {\it $i$th dual idempotent of $\Gamma$ with
respect to $x.$} For convenience we define $E_i^*=0$ for $i<0$ and
$i>D$. We observe
(i) $\sum_{i=0}^DE_i^*=I$; 
(ii) $\overline{E_i^*}=E_i^*$ $(0\leq i\leq D)$,
 (iii) ${E_i^*}^t=E_i^*$ $(0\leq i\leq D)$,
(iv) $E_i^*E_j^*=\delta_{ij}E_i^*$ $(0\leq i,j\leq D)$.
The $E^*_i$ have the following interpretation.
 Using
(\ref{29}) we find
\begin{equation}
\nonumber 
E_i^*V={\rm span}\{\hat y|y\in X,\ \ \partial (x,y)=i\}\ \ \ \ \
  (0\leq i\leq D).
\end{equation}
By this 
and since $\lbrace {\hat y} |y \in X\rbrace $ is an orthonormal
basis for $V$,
\begin{equation}
\nonumber 
  V=E_0^*V+E_1^*V+\cdots+E_D^*V\ \ \ \ \ ({\rm orthogonal\ direct\
  sum)}.
\end{equation}
For $0 \leq i \leq D$, $E^*_i$ acts on $V$ as the projection
onto $E^*_iV$. We call
$E_i^*V$  the {\it $i$th subconstituent of $\Gamma$
with respect to $x$.}
For $0 \leq i \leq D$ we 
define $s_i =\sum  {\hat y}$, where the sum is over all
vertices $y \in X$ such that $\partial(x,y)=i$.
We observe $s_i \in E^*_iV$.

\medskip
\noindent 
Let $\mathbf{T}=\mathbf{T}(x)$ denote the subalgebra of ${\rm
Mat}_X(\mathbb{C})$ generated by $A,$ $E_0^*,$ $E_1^*$, $\cdots$,
$E_D^*.$ The algebra 
$\mathbf{T}$ is semisimple but not commutative in general 
\cite[Lemma 3.4]{te92}.
We call $\mathbf{T}$ the {\it Terwilliger algebra} (or {\it subconstituent
algebra}) of
$\Gamma$ with respect to $x$.
  We refer the reader to
\cite{balm,
banmun,  
cau,
curtin99-1,
curtin99-2,
curtin4,
cnom,
egge1,
egge2,
GoTe,
hobart,
ishi1, 
ishi2,
sag,
tanabe,
te92,
te93,
te93-2,
terendone,
ternewineq,
tomyam}
for more
information on the Terwilliger algebra.
We will use the following facts.
Pick any integers $h,i,j$ $(0 \leq h,i,j\leq D)$.
By \cite[Lemma 3.2]{te92} we have
  $E_i^*A_hE_j^*=0$ if and only if
  $p_{ij}^h=0$.
By this and 
(\ref{2}),
(\ref{3_30.3}) we find
\begin{eqnarray}\label{34}
  E_i^*A_hE_1^*&=&0\ \ \ \ {\rm if}\ \ \ \ |h-i|>1
  \ \ \ \ \ \ (0\leq h,i\leq D),
  \\
  E_i^*AE_j^*&=&0\ \ \ \ {\rm if}\ \ \ \ |i-j|>1
  \ \ \ \ \ \ (0\leq i,j\leq D).
\label{3_30.1}
\end{eqnarray}

\begin{lem}\label{l0}
The following (i), (ii) hold for $0 \leq i\leq D$.
\begin{enumerate}
\item[(i)] $E_i^*JE_1^*=E_i^*A_{i-1}E_1^*+E_i^*A_iE_1^*+E_i^*A_{i+1}E_1^*.$
\item[(ii)]
$A_iE_1^*=E_{i-1}^*A_iE_1^*+E_i^*A_iE_1^*+E_{i+1}^*A_iE_1^*.$
\end{enumerate}
\end{lem}
\begin{proof} (i) 
Recall $J=\sum_{h=0}^D A_h$ so  $E^*_iJE^*_1
  =\sum_{h=0}^D E_i^*A_hE_1^*.$
Evaluating this using (\ref{34}) we obtain the result.
%
%
\\
\noindent (ii) Recall $I=\sum_{h=0}^DE^*_h$ so
$A_iE_1^*=\sum_{h=0}^DE_h^*A_iE_1^*.$
Evaluating this using (\ref{34}) we obtain the result.
\end{proof}


\begin{lem} For $0\leq i\leq D-1$ we have
\begin{equation}\label{l02e}
E_{i+1}^*A_iE_1^*-E_i^*A_{i+1}E_1^*=\sum\limits_{h=0}^iA_h E_1^*
- \sum\limits_{h=0}^iE_h^* JE_1^*.
\end{equation}
\end{lem}
\begin{proof}
%
Evaluate each term in 
the right-hand side of
(\ref{l02e}) 
using Lemma~\ref{l0}
and simplify the result.
\end{proof}

\begin{cor}
Let $v$ denote a vector in $E^*_1V$ which is orthogonal to $s_1$.
Then  for $0 \leq i\leq D-1$ we have
\begin{equation}
E_{i+1}^*A_iv-E_i^*A_{i+1}v=\sum\limits_{h=0}^iA_hv.
\label{l02en}
\end{equation}
\end{cor}
\begin{proof}
Apply all terms of (\ref{l02e})  to $v$ and evaluate the result
using $E^*_1v=v$  and
$Jv=0$.
\end{proof}

\begin{lem}\label{10-3}
The following (i), (ii) hold for $1 \leq i \leq D-1$.
\begin{enumerate}
\item[(i)]  $E_{i+1}^*AE^*_iA_{i-1}E^*_1=c_i E^*_{i+1}A_iE^*_1 $
\item[(ii)] $E_{i-1}^*AE^*_iA_{i+1}E^*_1=
b_iE^*_{i-1}A_iE^*_1.$
\end{enumerate}
\end{lem}
\begin{proof}(i)
For all $y,z\in X$, on either side the $yz$ entry is
equal to $c_i$ if $\partial(x,y)=i+1$,
 $\partial(x,z)=1$,
 $\partial(y,z)=i$, and zero otherwise.
\\
\noindent (ii)
For all $y,z\in X$, on either side the $yz$ entry is
equal to $b_i$ if $\partial(x,y)=i-1$,
 $\partial(x,z)=1$,
 $\partial(y,z)=i$, and zero otherwise.
\end{proof}

\begin{cor}\label{ll1}
Let $v$ denote a vector in $E^*_1V$.
 Then the following (i), (ii) hold for $1 \leq i \leq D-1$.
\begin{enumerate}
\item[(i)]
Suppose $E^*_iA_{i-1}v=0$. Then $E^*_{i+1}A_iv=0$.
\item[(ii)] 
Suppose $E^*_{i}A_{i+1}v=0$. Then $E^*_{i-1}A_{i}v=0$.
\end{enumerate}
\end{cor}
\begin{proof} In
Lemma~\ref{10-3}(i),(ii)  apply both sides to $v$ and 
use $E^*_1v=v$.
\end{proof}

\section{The modules of the Terwilliger algebra}\label{s5}

\noindent
Let $\mathbf{T}$ denote the Terwilliger algebra of $\Gamma$
with respect to
$x.$  
By a {\it $\mathbf{T}$-module} we mean a subspace $W\subseteq V$ such
that $BW\subseteq W$ for all $B\in \mathbf{T}.$
Let $W$ denote a $\mathbf{T}$-module. 
Then $W$ is said to be {\it irreducible} whenever
$W$ is nonzero and $W$ contains no $\mathbf{T}$-modules other than
$0$ and $W$.
Let $W$ denote an irreducible $\mathbf{T}$-module. 
Then 
$W$ is the orthogonal
direct sum of the nonzero spaces  among
$E^*_0W, E^*_1W, \ldots, E^*_DW$
\cite[Lemma 3.4]{te92}.
By the {\it endpoint} of
$W$ we mean 
  ${\rm min}\{i|0\leq i\leq D, E_i^*W\not=0\}$.
By the {\it diameter} of $W$ we mean 
  $|\{i|0\leq i\leq D, E_i^*W\not=0\}|-1$.
We say $W$ is {\it thin} whenever $E^*_iW$ has dimension at most 1
for 
$0\leq i\leq D$.
There exists a unique irreducible 
 $\mathbf{T}$-module  which has endpoint 0
\cite[Prop. 8.4]{egge1}. 
This module is called $V_0$.
For $0 \leq i \leq D$ the vector 
$s_i$ is a basis for 
$E^*_iV_0$
\cite[Lemma 3.6]{te92}.
Therefore 
$V_0$ is
thin with diameter $D$.
The module $V_0$ is orthogonal to each irreducible 
 $\mathbf{T}$-module other than $V_0$
\cite[Lem. 3.3]{curtin99-1}.
For more information
on $V_0$ see
\cite{curtin99-1, egge1}.
We will use the following facts.


\begin{lem}\label{l2.1}\cite[Lemma 3.9]{te92}
Let $W$ denote an irreducible $\mathbf{T}$-module
with endpoint $r$ and diameter $d.$ Then
\begin{equation}\label{41}
  E_i^*W\not=0\ \ \ \ (r\leq i\leq r+d).
\end{equation} Moreover
\begin{equation}\label{3-30.5}
E^*_iAE^*_jW\not=0\ \ \ \ {\rm if}\ |i-j|=1,\ \ (r\leq i, j\leq
r+d).
\end{equation}
\end{lem}
\begin{lem}\label{extra} \cite[Lemma~3.4]{curtin99-1} 
Let $W$ denote a $\mathbf{T}$-module.
Suppose there exists an integer $i$ $(0\leq i\leq D)$ such that
${\rm dim}(E_i^*W)=1$ and $W=\mathbf{T}E_i^*W$. Then $W$ is
irreducible.
\end{lem}
\begin{thm}\label{t3-4-6}\cite[Lemma~10.1]{GoTe}, 
\cite[Theorem 11.1]{terendone}
Let $W$ denote a thin irreducible $\mathbf{T}$-module with
endpoint one, and let $v$ denote a nonzero vector in $E_1^*W.$
Then $W=\mathbf{M}v.$ Moreover the diameter of $W$ is $D-2$ or $D-1$.
\end{thm}
\begin{thm}\label{t4.5.1} \cite[Corollary 8.6,  Theorem 9.8]{GoTe}
Let $v$ denote a nonzero vector in $E_1^*V$ which is
orthogonal to
$s_1$.
Then the dimension of
$\mathbf{M}v$ is $D-1$ or $D$. 
 Suppose
the dimension of $
 \mathbf{M}v$ is $D-1$. Then
$\mathbf{M}v$ is a thin irreducible $\mathbf{T}$-module with
endpoint $1$ and diameter $D-2.$
\end{thm}
%
%
%
%

\section{The proof of Theorem~\ref{10-20i}}\label{main1}

\noindent We now give a proof of
Theorem 1.1.

\begin{proof}
((i) $\Longrightarrow$ (ii)) 
We show
$\mathbf{M}v$ 
is a thin irreducible 
 $\mathbf{T}$-module with
 endpoint 1.
By 
Theorem \ref{t4.5.1} 
 the dimension of
$\mathbf{M}v$ is either $D-1$ or $D$. First assume the dimension
of $\mathbf{M}v$ is equal to $D-1$. Then by Theorem~\ref{t4.5.1},
$\mathbf{M}v$ is a thin irreducible $\mathbf{T}$-module with
endpoint 1. Next assume the dimension of
$\mathbf{M}v$ is equal to $D$. The space $(\mathbf{M};v)$ contains $J$ and
has dimension at least 2, so there exists $P \in (\mathbf{M};v)$
such that $J, P$ are linearly independent. From the construction
$Pv \in E^*_DV$. Observe $Pv \not=0$; otherwise the dimension of
$\mathbf{M}v$ is not $D$. The elements $A_0, A_1, \ldots, A_D$
form a basis for $\mathbf{M}$. Therefore the elements
$A_0+A_1+\cdots + A_i$ $ (0 \leq i \leq D)$ form
 a basis for $\mathbf{M}$.
Apparently there exist complex scalars $\rho_i$ $(0 \leq i \leq D)$ such
that $P=\sum_{i=0}^D \rho_i (A_0+A_1 +\cdots + A_i)$. Recall
$J=\sum_{h=0}^D A_h$. Subtracting a scalar multiple of $J$
from $P$ if necessary, we may assume $\rho_D=0$. We consider $Pv$
from two points of view.
 On one hand
we have $Pv \in E^*_DV$. Therefore $E^*_DPv = Pv $ and
$E^*_iPv=0$ for $0 \leq i \leq D-1$. On the other hand
using (\ref{l02en}),
$$
Pv = \sum_{i=0}^{D-1} \rho_i (E^*_{i+1}A_iv-E^*_iA_{i+1}v).
$$
Combining these two points of view we find
$Pv=\rho_{D-1}E^*_DA_{D-1}v$, $E^*_0Av=0$,   and
\begin{eqnarray}
\rho_{i-1} E^*_iA_{i-1}v = \rho_iE^*_iA_{i+1}v \qquad (1 \leq i
\leq D-1). \label{eq:rhodep}
\end{eqnarray}
We mentioned $Pv\not=0$; therefore
 $\rho_{D-1}\not=0$ and
 $E^*_DA_{D-1}v \not=0$.
Applying Corollary \ref{ll1}(i) we find $E^*_iA_{i-1}v \not=0$ for $1
\leq i \leq D$. We claim $E^*_iA_{i+1}v$ and $E^*_iA_{i-1}v$ are
linearly dependent for $1 \leq i \leq D-1$. Suppose there exists
an integer $i$ $(1 \leq i \leq D-1)$ such that
 $E^*_iA_{i+1}v$
and $E^*_iA_{i-1}v$ are linearly independent. Then
 $E^*_iA_{i+1}v \not=0$. Applying
Corollary
\ref{ll1}(ii) we find
 $E^*_jA_{j+1}v \not=0$
for $i \leq j \leq D-1$. Using these facts and
(\ref{eq:rhodep}) we routinely find
$\rho_j=0$ for $i \leq j \leq D-1$. In particular $\rho_{D-1}=0$
for a contradiction. We have now shown $E^*_iA_{i+1}v$ and
$E^*_iA_{i-1}v$ are linearly dependent for $1 \leq i \leq D-1$.
Observe $\mathbf{M}v$ is spanned by the vectors
$$
(A_0+A_1+\cdots + A_i)v \qquad \qquad (0 \leq i \leq D-1).
$$
By (\ref{l02en}) and our above
comments we find $\mathbf{M}v$ is contained in the span of
\begin{eqnarray}
E^*_{i+1}A_iv \qquad \qquad (0 \leq i \leq D-1). \label{eq:span}
\end{eqnarray}
Since $\mathbf{M}v$ has dimension $D$ we find $\mathbf{M}v$ is
equal to the span of (\ref{eq:span}). Apparently $\mathbf{M}v$ is a
$\mathbf{T}$-module. Moreover $\mathbf{M}v$ is irreducible by
Lemma~\ref{extra}. Apparently  $\mathbf{M}v$ is thin with endpoint
1.
\medskip

\noindent ((ii) $\Longrightarrow$ (i))
We show 
 $(\mathbf{M};v)$ has 
dimension at least 2.
Since 
 $J \in (\mathbf{M};v)$ it suffices to exhibit an element
 $P \in (\mathbf{M};v)$ such that $J, P$ are linearly independent.
Let $W$ denote a thin
irreducible $\mathbf{T}$-module which has endpoint 1 and contains
$v$. By Theorem~\ref{t3-4-6} we have $W=\mathbf{M}v$; also by
Theorem~\ref{t3-4-6} 
the diameter of $W$ is $D-2$ or $D-1$. First
suppose $W$ has diameter $D-2$. 
Then $W$ has dimension $D-1$.
Consider the map $\sigma
:\mathbf{M} \rightarrow V$ which sends each element $P$ 
to $Pv$. The image of 
$\mathbf{M}$ 
under 
$\sigma $ is $\mathbf{M}v$ and
the kernel of $\sigma$ is contained in $(\mathbf{M};v)$. The image
has dimension $D-1$ and $\mathbf{M}$ has dimension $D+1$ so the
kernel has dimension $2$. It follows $(\mathbf{M};v)$ has
dimension at least 2. Next assume $W$ has diameter $D-1$. In this
case $E^*_DW\not=0$ by (\ref{41}). Since $W=\mathbf{M}v$ there
exists $P \in \mathbf{M}$ such that $Pv $ is a nonzero element in
$E^*_DW$. Now $P \in (\mathbf{M};v)$. Observe $P,J$ are linearly
independent since $Pv\not=0$ and $Jv=0$.
Apparently the dimension of $(\mathbf{M};v)$ is
at least 2.

\medskip
\noindent  Now assume (i), (ii) hold. We show the dimension of
$(\mathbf{M};v)$ is 2. To do this, we show the dimension of
$(\mathbf{M};v)$ is at most 2.
Let $H$ denote the subspace of
$\mathbf{M}$ spanned by $A_0, A_1,\ldots, A_{D-2}$. 
We show $H$ has
 0 intersection with $(\mathbf{M};v)$.
By Theorem
\ref{t4.5.1} 
the dimension of $\mathbf{M}v$ is at least $D-1$.
Recall
$\mathbf{M}$ is generated by $A$ so
the vectors $A^iv$ $(0 \leq i \leq D-2)$ are linearly independent.
Apparently the vectors 
 $A_iv$ $(0 \leq i \leq D-2)$ are linearly independent.
For $0 \leq i \leq D-2$ the vector $A_iv$
is contained  in $\sum_{h=0}^{D-1}E^*_hV$
by 
Lemma~\ref{l0}(ii); therefore $A_iv$ is orthogonal to $E^*_DV$.
We now see the vectors   
 $A_iv$ $(0 \leq i \leq D-2)$ are linearly independent
 and orthogonal to $E^*_DV$. It follows
 $H$ has 0 intersection with $(\mathbf{M};v)$.
Observe $H$ is codimension 2 in
$\mathbf{M}$ so the dimension of $(\mathbf{M};v)$ is at most 2. We
conclude the dimension of
 $(\mathbf{M};v)$ is 2.
\end{proof}

\section{Pseudo primitive idempotents}\label{s3}

\noindent In this section we introduce the notion of
a pseudo primitive
idempotent.
\begin{defn}
\label{def:ppi} 
 For each $\theta\in \mathbb{C} \cup
\infty$ we define a  subspace
of
 $\mathbf{M}$ which we call
$\mathbf{M}(\theta)$.
For $\theta\in \mathbb{C}$, $\mathbf{M}(\theta)$ consists of those
elements $Y $ of $\mathbf{M} $ such that $(A-\theta I)Y\in
\mathbb{C} A_D$.
We define $\mathbf{M}(\infty)=
\mathbb{C}A_D.$
\end{defn}
\noindent With reference to 
Definition
\ref{def:ppi}, 
we will show
each $\mathbf{M}(\theta)$ has dimension 1.
To establish this we display a basis for
$\mathbf{M}(\theta)$.
We will use the following result.
\begin{lem} 
\label{lem:rhorec}
Let $Y$ denote an element of
$\mathbf{M}$ and write
$Y=
\sum\limits_{i=0}^D\rho_iA_i.$ Let $\theta$ denote a complex
number. Then the following (i), (ii)  are equivalent.
\begin{enumerate}
\item[(i)]  $(A-\theta I)Y \in
           \mathbb{C}A_D.$
\item[(ii)]
$\rho_i = \rho_0 f_i(\theta)k^{-1}_i $ for $0 \leq i\leq D$.
\end{enumerate}
\end{lem}
\begin{proof} Evaluating $(A-\theta I)Y$
 using
$Y=\sum\limits_{i=0}^D\rho_iA_i$ and  simplifying the result using
(\ref{10}) we obtain
$$
(A-\theta I)Y =
\sum\limits_{i=0}^DA_i(c_i \rho_{i-1}+ a_i \rho_i + b_i \rho_{i+1} - \theta \rho_i),
$$
where $\rho_{-1}=0$ and $\rho_{D+1}=0$.
Observe by (\ref{4}), (\ref{14}) that $\rho_i = \rho_0
f_i(\theta)k^{-1}_i$ for $0 \leq i\leq D$ if and only if $c_i
\rho_{i-1}+ a_i \rho_i + b_i \rho_{i+1} = \theta \rho_i$ for $0
\leq i \leq D-1$. The result follows.
\end{proof}

\begin{cor}
\label{cor:basism15}
For $\theta \in
 \mathbb{C}$ the following is a basis for
$\mathbf{M}(\theta)$.
\begin{equation}\label{e5-27}
\sum\limits_{i=0}^Df_i(\theta)k^{-1}_iA_i.
\end{equation}
\end{cor}
\begin{proof} Immediate from 
Lemma
\ref{lem:rhorec}.
\end{proof}

\begin{cor}
\label{cor:dim1}
The space
$\mathbf{M}(\theta)$ has dimension 1
for all $\theta\in \mathbb{C}\cup \infty$.
\end{cor}
\begin{proof}
Suppose
 $\theta =\infty$. Then
$\mathbf{M}(\theta)$ has basis $A_D$ and therefore  has dimension 1.
Suppose
$\theta\in \mathbb{C}.
$
Then
$\mathbf{M}(\theta)$ has  dimension 1 by
Corollary \ref{cor:basism15}.
\end{proof}

\begin{lem}\label{l5-20} Let $\theta$ and
$\theta'$ denote distinct elements of
$\mathbb{C}\cup \infty$. Then
$\mathbf{M}(\theta) \cap
\mathbf{M}(\theta')=0$.
\end{lem}
\begin{proof} 
This is a routine consequence of
Corollary \ref{cor:basism15} and the fact that $
\mathbf{M}(\infty)=\mathbb{C} A_D$.
\end{proof}

\begin{cor}\label{cor:add}
For $0 \leq i \leq D$ we have
$\mathbf{M}(\theta_i)=\mathbb{C}E_i$.
\end{cor}
\begin{proof}
Observe 
$(A-\theta_iI)E_i=0$ 
so $E_i \in
\mathbf{M}(\theta_i)$.
The space 
$\mathbf{M}(\theta_i)$ has dimension 1 by
Corollary
\ref{cor:dim1}
and 
$E_i$ is nonzero
so $E_i$ is a basis for
$\mathbf{M}(\theta_i)$.
\end{proof}
\begin{remark}
\cite[p. 63]{BanIto}
For $0 \leq j \leq D$
we have
$$
E_j = m_j |X|^{-1}
\sum\limits_{i=0}^D f_i(\theta_j)k^{-1}_iA_i,
$$
where $m_j$ denotes the rank of $E_j$. 
\end{remark}
\begin{defn}
Let $\theta \in 
\mathbb{C}\cup \infty$.
By a {\it pseudo primitive idempotent} for $\theta$ we mean
a nonzero element of
$\mathbf{M}(\theta)$, where 
$\mathbf{M}(\theta)$ is from
Definition
\ref{def:ppi}.
\end{defn}

\section{The local eigenvalues}\label{s6}

\begin{defn}\label{10-3-2} Define a function\ \
$\widetilde{\
}:\mathbb{C}\cup \infty\longrightarrow\mathbb{C}\cup\infty$
by $$\widetilde{\eta}=
 \begin{cases}
    \infty\ \ \  & \text{if $\eta=-1$},\\
    -1\ \ \      & \text{if $\eta=\infty$},\\
    -1-\frac{b_1}{1+\eta}\ \ \  & \text{if $\eta \not=-1, \eta\not=\infty$}.
  \end{cases}
$$ Observe $\widetilde{\widetilde{\eta}}=\eta$\
for all $\eta\in\mathbb{C}\cup \infty.$
\end{defn}
\noindent 
Let $v$ denote a nonzero vector in $E^*_1V$ which is orthogonal to
$s_1$. Assume $v$ is an eigenvector for $E^*_1AE^*_1$
and let 
$\eta $ denote the corresponding 
eigenvalue.
 We recall a few facts  concerning
$\eta $ and 
${\widetilde \eta}$.
We have ${\widetilde \theta_1} \leq
\eta \leq {\widetilde \theta_D}$
 \cite[Theorem 1]{te86}.
If $\eta=
{\widetilde \theta_1}$ then
 $\widetilde \eta =\theta_1$.
If $\eta= {\widetilde \theta_D}$ then
 $\widetilde \eta =\theta_D$.
We have
$\theta_D < -1 < \theta_1$ by 
 \cite[Lemma 3]{te86} so
${\widetilde \theta_1} < -1
< {\widetilde \theta_D}$.
If ${\widetilde \theta_1} < \eta < -1 $ then 
$\theta_1 < {\widetilde \eta}.$
 If 
$-1 < \eta < {\widetilde \theta_D} $
then
${\widetilde \eta} < \theta_D$. We will show that if
${\widetilde \theta_1} < \eta <
 {\widetilde \theta_D} $
 then
${\widetilde \eta}$ is not an eigenvalue of $\Gamma$. Given the
above inequalities, to prove this it suffices to prove the
following result.
\begin{prop}\label{p5}
Let $v$ denote a nonzero vector in $E_1^*V.$ 
Assume $v$ is an eigenvector
for $E^*_1AE^*_1$ and let $\eta$ denote the corresponding
eigenvalue. Then $\widetilde{\eta}\not=k$.
\end{prop}
\begin{proof} Suppose $ \widetilde{\eta}=k.$ Then 
$\eta = 
 \widetilde{k}$ so by 
Definition~\ref{10-3-2}, 
$$\eta=-1-\frac{b_1}{k+1}.$$
By this and since $b_1<k$ we see
$\eta$ is a rational number such that
$-2<\eta < -1$.
In particular $\eta$ is not an integer.
Observe $\eta$ is an
eigenvalue of the subgraph of $\Gamma$
induced on the set of vertices adjacent
$x;$ therefore $\eta$ is an algebraic integer.
A rational algebraic integer is an
integer so we have a contradiction.
We conclude 
$\widetilde{\eta}\not=k$.
\end{proof}

\begin{cor}\label{5-27}
Let $v$ denote a nonzero vector in $E_1^*V$ which is orthogonal to $s_1$.
Assume $v$ is an eigenvector
for $E^*_1AE^*_1$ and let $\eta$ denote the corresponding
eigenvalue. Suppose ${\widetilde \theta_1} < \eta <
 {\widetilde \theta_D} $.
 Then
${\widetilde \eta}$ is not an eigenvalue of $\Gamma$.
\end{cor}
%
%
%
%

\section{The proof of Theorem~\ref{11-23}}\label{main2}
\bigskip

\noindent 
We now give a  proof of
Theorem~\ref{11-23}.

\begin{proof}
We first show $E$ is contained in $(\mathbf{M};v)$. To do this we
show $Ev \in E^*_DV$. 
First
suppose $\eta\not=-1.$ Then 
$\widetilde{\eta}
\in\mathbb{C}$ by
Definition~\ref{10-3-2}.
By Definition
\ref{def:ppi}
there exists $\epsilon\in \mathbb{C}$ such that
$(A-{\widetilde {\eta}} I)E=\epsilon
A_D.$ By this and 
 Lemma~\ref{l0}(ii),
\begin{eqnarray}\label{3_30}
AEv&=&{\widetilde \eta} Ev+\epsilon A_Dv \nonumber\\
   &\in& \mathbb{C}Ev+E_{D-1}^*W+E_D^*W.
\end{eqnarray}
In order to show 
$Ev \in E^*_DV$ 
we show
$E_i^*Ev=0$ for $0\leq i\leq D-1$.
Observe $E^*_0Ev=0$ since
$E^*_0Ev \in E^*_0W$ and $W$ has endpoint 1. We show $E^*_1Ev=0.$
By Corollary
\ref{cor:basism15}
there exists a nonzero $m 
\in\mathbb{C}$ 
 such that
$$E=m\sum_{h=0}^Df_h({\widetilde \eta})k^{-1}_hA_h.$$
 Let us
abbreviate
\begin{equation}
\label{eq:extra}
\rho_h = mf_h({\widetilde \eta})k^{-1}_h \qquad (0 \leq h \leq D),
\end{equation}
so that $E=\sum_{h=0}^D \rho_hA_h$.
By this and
(\ref{34})
we find
 $E^*_1EE^*_1=\sum_{h=0}^2 \rho_hE^*_1A_hE^*_1$. Applying this to $v$
 we find
\begin{equation}\label{p1}
E_1^*Ev=\sum_{h=0}^2 \rho_h E^*_1A_hv.
\end{equation}
Setting $i=1$ in
Lemma~\ref{l0}(i), applying each term to $v$, and using $Jv=0$ we find
\begin{equation}\label{p2}
0=\sum_{h=0}^2 E^*_1A_hv.
\end{equation}
By (\ref{p1}), (\ref{p2}), and since $E^*_1Av=\eta v$ we find
$E^*_1Ev=\gamma v$ where 
$\gamma = \rho_0-\rho_2+\eta(\rho_1-\rho_2)$.
Evaluating $\gamma$ 
using
(\ref{14}), 
(\ref{eq:extra}), 
 and Definition \ref{10-3-2} we routinely find $\gamma=0$.
Apparently
$E_1^*Ev=0$.
 We now show $E^*_iEv=0$ for $2\leq i\leq D-1.$
Suppose there exists an integer
$j$ $(2 \leq j \leq D-1)$
such that 
$E^*_jEv\not=0$. We choose $j$ minimal so that
\begin{equation}\label{4-4-1}
E_i^*Ev=0\ \ \ \ (0\leq i\leq j-1).
\end{equation}
Combining this with
(\ref{3_30}) we find
\begin{equation}\label{3_30.2}
E_i^*AEv=0\ \ \ \ (0\leq i\leq j-1).
\end{equation}
Since $W$ is thin and since $E^*_jEv\not=0$ we find
$E^*_jEv$ is a basis for
$E^*_jW.$ 
Apparently $E^*_{j-1}AE^*_jEv$ spans $E^*_{j-1}AE^*_jW$.
The space
$E^*_{j-1}AE^*_jW$ is
 nonzero by (\ref{3-30.5}) and since the diameter of $W$ is at least $D-2$.
Therefore $E^*_{j-1}AE^*_jEv\not=0$.
We may now argue
$$\begin{array}{ccll} E^*_{j-1}AEv
&=&\sum_{i=0}^DE^*_{j-1}AE^*_iEv & \cr
&=& E^*_{j-1}AE^*_jEv  &{\rm by}\ (\ref{3_30.1}), (\ref{4-4-1})\cr
\ &\not=&0
\end{array}$$
which contradicts (\ref{3_30.2}).
We conclude $E^*_iEv=0$ for $2 \leq i \leq D-1$.
We have now shown
 $E^*_iEv=0$ for $0 \leq i \leq D-1$ so
$Ev \in E^*_DV$ in the case $\eta \not=-1$.
Next suppose $\eta=-1$, so that ${\widetilde \eta}=\infty$. By
Definition
\ref{def:ppi} there exists a nonzero  $t \in {\mathbb C}$  such that 
$E=tA_D$.
In order to show $Ev \in E^*_DV$ we show
$A_Dv
\in E^*_DV$.
Since $A_Dv$ is contained in $E^*_{D-1}V+E^*_DV$
by Lemma~\ref{l0}(ii), it suffices to  show
$E^*_{D-1}A_Dv=0.$
 To do this it is convenient to prove a bit more,
that $E^*_iA^*_{i+1}v=0$ for $1\leq i\leq D-1.$ We prove this by
induction on $i.$ First assume $i=1$.
Setting $i=1$ in Lemma~\ref{l0}(i), applying
each term to $v$ and using $Jv=0$,  $E^*_1Av=-v,$ we obtain 
$E^*_1A^*_2v=0$.
 Next suppose $i \geq 2$ and assume by induction that
$E^*_{i-1}A_iv=0$. We show 
$E^*_iA_{i+1}v=0.$
To do this we assume 
$E^*_iA_{i+1}v\not=0$ and get a contradiction.
 Note that $E^*_iA_{i+1}v$ spans $E_i^*W$
since $W$ is thin. Then $E^*_{i-1}AE^*_iA_{i+1}v\not=0$ by
(\ref{3-30.5}). But $E^*_{i-1}AE^*_iA_{i+1}v=b_iE^*_{i-1}A_iv$ by
Lemma~\ref{10-3}(ii). Of course $b_i\not=0$ so
$E^*_{i-1}A_iv\not=0,$ a contradiction. 
Therefore 
$E^*_iA^*_{i+1}v=0$.
We have now shown
$E^*_iA^*_{i+1}v=0$ for $1\leq i\leq D-1$ and in particular
$E^*_{D-1}A_Dv=0.$
It follows 
$Ev \in E^*_DV$ for the case $\eta =-1$.
We have now shown
$Ev \in E^*_DV$ for all cases
so
$E\in (\mathbf{M};v)$. We now prove 
 $E,J$ form a basis for $(\mathbf{M};v)$. By Theorem 1.1
$(\mathbf{M};v)$ has dimension 2. We mentioned earlier $J \in
(\mathbf{M};v)$. We show $E,J$ are linearly independent. 
Recall $E, J$ are pseudo primitive idempotents
for $\widetilde
\eta, k$ respectively.
We have
${\widetilde{\eta}}\not=k$ by Proposition~\ref{p5}
so $E,J$ are linearly independent in view of Lemma~\ref{l5-20}.
\end{proof}

\noindent {\bf Acknowledgements} The initial work for this paper
was done when the second author was an Honorary Fellow at the
University of Wisconsin-Madison (July-December 2000) supported by
the National Science Council, Taiwan ROC.

{\small

\noindent Paul Terwilliger \\
 Department of Mathematics \\ 
University of Wisconsin \\
480 Lincoln Drive \\
Madison Wisconsin \\
USA  53706 \\
Email: terwilli@math.wisc.edu\\

\medskip
\noindent Chih-wen Weng\\
Department of Applied 
Mathematics\\
National Chiao Tung University\\
1001 Ta Hsueh Road\\
Hsinchu\\
Taiwan ROC\\
Email: weng@math.nctu.edu.tw

\end{document}